\documentclass[11pt,a4paper,leqno]{article}
\usepackage{latexsym}
\usepackage{amsfonts}
\usepackage {amsbsy}
\usepackage {amsmath}
\usepackage{amssymb}
\newtheorem{defn}{Definition}[section] 
\newtheorem{thm}[defn]{Theorem} 
\newtheorem{prop}[defn]{Proposition} 
\newtheorem{lem}[defn]{Lemma} 
\newtheorem{cor}[defn]{Corollary}

\newcommand \psh{plurisubharmonic }
\newcommand \demo{ Proof: }
\newcommand \C{\mathbb C}
\newcommand \N{\mathbb N}

\newcommand \R{\mathbb R}
\newcommand \B{\mathbb B}
\newcommand \fin{$\blacktriangleright$\\}
\newcommand \lto{\longmapsto}
\newcommand \lra{\longmapsto}

\newcommand \mc{\mathcal}
\newcommand \mb{\mathbb}
\newcommand \mrm{\mathrm}
\newcommand \W{\Omega}
\newcommand \w{\omega}
\newcommand \sm{\setminus}
\newcommand \sub{\subset}
\newcommand \Sub{\Subset}
\newcommand \al{\alpha}

\newcommand \de{\delta}

\newcommand \Fi{\Phi}

\newcommand \la{\lambda}

\newcommand \ka{\kappa}
\newcommand \vep{\varepsilon}
\newcommand \vfi{\varphi}
\newcommand \si{\sigma}

\newcommand  \vt{\vartheta}
\newcommand \ti {\Tilde}

\newcommand \ove{\overline}
\newcommand \bd{\partial}
\newcommand \we{\wedge}

\numberwithin{equation}{section}
 \title{ A minimum principle \\ for 
 plurisubharmonic functions}
\author{Ahmed  Zeriahi} 
\date{}
\begin{document}
\maketitle

\noindent {\bf Abstract: } The main goal of this work is to give  new and precise generalizations to various classes of plurisubharmonic functions of the classical minimum modulus principle for holomorphic functions of one complex variable, in the spirit of the famous lemma of Cartan-Boutroux. As an application we obtain precise estimates on the size of "plurisubharmonic lemniscates" in terms of appropriate Hausdorff contents. 
\vskip 0.3 cm 
 \noindent {\bf AMS 2000 Mathematical subject classification:} 31C10, 31C15,  
 32F05, 32F99, 32U05, 32U99.
 \vskip 0.3 cm 
\section {Introduction} 

  In Complex Analysis of one variable, besides the classical 
 maximum principle, there is another principle that is less well 
 known, but just as important and somewhat more subtil. It consists 
 in giving a precise lower bound for the modulus of a holomorphic 
 function (normalized in a suitable way) on a given open disc at 
 all points of a smaller disc, except those which belong to an 
 exceptionnal subset containing its zeros, in terms of its maximum 
 on the given disc; the size of the exceptionnal set should be 
 precisely estimated in terms of capacity or one dimensionnal  
 Hausdorff content. This is called the minimum modulus principle  
 for holomorphic functions (see [26], [35]).
  
  This principle plays an important role in many problems involving  rationnal or meromorphic functions which may have many poles in a   given domain and then it's desirable to find an upper bound of 
 such function, which remains to control precisely the subset where the  
 denominator of the function is small.
 
 Classically in one complex variable, this can be done using the 
 well know Cartan-Boutroux' lemma, which gives a precise uniform 
 estimate of the size of (monic) polynomial lemniscates in terms of 
 one dimensionnal Hausdorff content.
 
 This kind of estimates can be applied in Nevalinna theory 
 for the study of the growth of holomorphic functions (see [10],   
 [28]) as well as in Pad\'e approximation when
 studing the convergence of rationnal Pad\'e approximants of 
 meromorphic functions (see [13], [18]) and also in Harmonic  
 Analysis and PDE's (see [16], [9]).
 
 Uniform estimates on the size of sublevel sets of some classes of 
 plurisubharmonic functions, called {\it \psh lemniscates}, have 
 been obtained in our earlier papers (see [39],[40], [8]).
 
 The main purpose of the present paper is to give precise and new 
 various generalizations of the Cartan-Boutroux' lemma to  
 logarithmic potentials in $\C^n$ as well as a general version of  
 the minimum principle for arbitrary \psh functions on euclidean  
 balls in $\C^n$. 
 
  From the version of the minimum principle for logarithmic  
 potentials in $\C^n$, we deduce a new comparison inequality  
 between appropriate Hausdorff contents and 
 the classical logarithmic capacity in $\C^n$.
 
  The method used here, known as the method of "excluding balls"  ("boules d'exclusion" in french), seems to be quite classical in  
 real potential theory when estimating integral potentials but  
 surprinsingly, it has never been used in our context although it  
 has been used in a different context by Avanissian ([4]). 
  
  It turns out that the same method allows us to obtain a quite 
 general lower estimate for \psh functions on the unit ball, which 
 implies a sort of "three-circle minimum principle" for \psh
 functions, which can be seen as the dual counterpart of 
 Hadamard three-circle inequalities.  
 
 Actually, the same method can be also applied on a compact 
 K\"ahler manifold to give a minimum principle for 
 quasiplurisubharmonic functions on the manifold.
  
 \section{The Cartan-Boutroux lemma and the minimum principle} 
 Let us first recall the famous classical lemma of  Cartan-Boutroux 
 (see [10]).  Let $P (z)$ be a monic polynomial of one complex  
 variable of degree $d \geq 1$. For any given $\vep > 0,$  
 consider the polynomial $\vep-$lemniscate of $P$ defined by 
 $E (P;\vep) := \{z \in \C \ ; \ \vert P (z) \vert \leq \vep^d \}.$ 

 Then there exists a finite covering of $E (P;\vep)$ by $d$
 open discs with radii $(r_{j})_{1 \leq j \leq d}$ satisfying the  
 estimate 
 \begin{equation}
 \sum_{1 \leq j \leq d} r_{j} \leq 2 e \vep.
 \end{equation}
 In other words $\log \vert P (z) \vert \geq - d \log (1 \slash  
 \vep)$ for all $z \in \C$ outside the union of $d$ open discs with  radii $(r_{j})_{1 \leq j \leq d}$ satisfying the estimate 
 $\sum_{1 \leq j \leq d} r_{j} \leq 2 e \vep.$

 From this estimate it is possible to derive the following minimum  
 principle for holomorphic functions (see [26], [35]).
 If $f$ is a holomorphic function on the disc $\{ z \in \C ; \vert 
 z\vert \leq 2 e R\}$ such that $f (0) = 1$. Then for any real  
 number $0 < \eta < 1$, the following lower bound
 $$ \log \vert f (z)\vert > -  H (\eta) \log M_f (2 e R), \ \  
 \mrm{where} \ \ H (\eta) := \log \Bigl(3 e^3 \slash 2 \eta\Bigr),$$
 holds for $\vert z\vert \leq R$ outside the union of a finite  
 number discs of radii $(r_j)$ with $\sum_j r_j \leq 2 \eta R$.

 There is a more general version of Cartan-Boutroux' lemma which 
 can be stated as follows (see [26]).
 For any $0 < \alpha \leq 2$ there exists a finite covering of $E  
 (P;\vep)$ by $d$ open discs with radii $(r_{j})$ satisfying the 
 estimate 
 \begin{equation}
 \sum_{j = 1}^d r_{j}^\alpha \leq e (2 \vep)^{\al}.
 \end{equation}
 In other words this means that for any $\vep \in ]0,1]$ the  
 lower bound $\vert P (z)\vert \geq \vep^d$ holds for all $z$  
 ouside the union of $d$ exceptional discs with radii 
 $(r_j)$ satisfying the estimate $\sum_{j} r_{j}^\alpha \leq 
 e (2  \vep)^{\al}$. 

 This is equivalent to the following estimate in terms of the  
 Hausdorff  
 content of dimension $\al$
 $$h^{\al} (E (P;\vep)) \leq e (2 \vep)^\alpha, \forall \vep \in  
 ]0, 1].$$
 
 Let us end this section by recalling the definition of the 
 Hausdorff contents in a more general setting, since it will be 
 used  later. 
 
 Let $(X,d)$ be a metric space and $p > 0$ a real number. Then for 
 a given real number $\de > 0$, by definition, the 
 $\de-$Hausdorff content of dimension $p$ of a subset $E \sub X$ is 
 defined as follows
 $$h^p_{\de} (E) := \inf \{ \sum_{j \in \N} r (B_j)^p ; E \sub 
 \bigcup_{j \in \N} B_j, \ (B_j)_{j \in \N} \in \mc B_{\de} 
 (X,d)\},$$
 where $\mc B_{\de} (X,d)$ is the class of all countable coverings 
 $(B_j)_{j \in \N}$ of the set $E$ by balls of the metric space  
 $(X,d)$ of radii at most $\de$ and $r (B_j)$ is the radius of the  
 ball $B_j$ for each $j \in \N$.  
 
  Usually as in the Cartan-Boutroux lemma, we could just take $\de  
 =  + \infty$ to be infinite, which means that in the 
 above definition, we don't ask for an priori bound on the radii 
 $(r_j)$ of the balls $(B_j)$ which cover $E$. The corresponding 
 number is denoted by $h^p (E) = h^p_{\infty} (E)$ and 
 called the Hausdorff content of dimension $p$ of the set $E$.
 
  Observe that the $p-$dimensional Hausdorff measure of the set $E$  is defined by $H^p (E) := \sup_{\de > 0} h^p_{\de} (E) = \lim_{\de 
  \downarrow 0} h^p_{\de} (E).$
  
 \section{Projective masses and Lelong numbers} 
 
 Let us recall some well known definitions and properties 
 concerning Lelong numbers ([25], [15], [19]).
 
 Let $\W$ be a domain in $\C^n$ and $PSH (\W)$ be the cone of \psh  
 functions $u$ on $\W$ such that $u \not \equiv  - \infty.$ 
 Then $PSH (\W) \subset L_{loc}^{1} (\W)$ is a 
 closed subset for the $L_{loc}^{1}-$topology 
 and then it is a complete metric space (see [19]). 

  Let us consider the usual differential operators 
 on $\C^n$ defined by $d = \partial + \overline {\partial}$ and 
 define $d^c := ({1 \slash 2 \pi} i) (\partial - \overline  
 {\partial})$ so that $dd^c = ({i \slash \pi}) \partial \overline 
 {\partial}$.

 This normalization is choosen so that the following Monge-Amp\`ere  equation holds
 $$ (dd^c \log \vert z \vert)^n = \de (z)$$
 in the sense of currents on $\C^n$, where $\de (z)$ is the Dirac  
 point-mass at the origin.

 Now recall that if $V \in PSH (\W)$ then $dd^c V$ is a closed  
 positive current of bidegree (1,1) on $\W$ (see [25]).
 
 For any fixed $a \in \W$ and $0 < r << 1$ such that $\B (a;r) := 
 \{ z \in \C^n \ ; \ \vert z - a \vert \leq r \} \Sub \W$, define    the projective mass of the current $dd^c V$ on the ball $\B (a,r)$
 as follows
 \begin{equation}
        \vt_{V} (a,r) := \int_{\B (a,r)} dd^c V \wedge
        (dd^c \log \vert z - a \vert)^{n - 1}.
        \label{eq:pmasse}
 \end{equation}
 Then by a well know result of Lelong ([25]), the following formula 
\begin{equation}
 \vt_{V} (a,r) = \frac{(n - 1)!}{ \pi^{n - 1} r^{2 n - 2}} 
 \int_{\B (a,r)} dd^c V \wedge {\beta}_{n - 1} 
 = \frac{\mu_{V} (\B (a,r))}{ \tau_{2 n - 2} r^{2 n - 2}}
        \label{eq:flelong}
\end{equation}
 holds, where $\tau_{2 n - 2}$ is the $ (2 n - 2)-$dimensional volume of the euclidean unit ball in $\C^{n - 1}$, $\textstyle \beta := (i \slash 2) \partial  
 \overline{\partial} \ \vert z \vert^{2},$ ${\beta}_{n - 1} :=   
 {\beta}^{n - 1} \slash (n - 1)! $ and $\mu_{V} := (1 \slash  2  
 \pi) \Delta V$ is the Riesz measure associated to $V.$  
  
 Then the projective mass of the current $dd^c V$ at the point $a$  
 is defined by the following formula
 \begin{equation}
 \vt_{V} (a) := \lim_{r \to 0^{+}} \vt_{V} (a,r) 
 = \lim_{r \to 0} 
 \frac{\mu_{V} (\B (a,r))}{ \tau_{2 n - 2} r^{2 n - 2}}.
          \label{}
 \end{equation}
 The positive number $\vt_{V} (a)$ is also 
 called {\it the Lelong number} of the current $dd^c V$ or the 
 Lelong number of the function $V$ at the point $a$. 

 By a classical result of V. Avanissian (see [3], [20]), the  
 Lelong number can be also expressed by the following formulas 
 \begin{eqnarray}
 \vt_{V} (a) & = & \lim_{r \to 0^{+}} 
 \frac {1} {\log r} {\int_{\vert \xi \vert = r}  
 V (a + r \  \xi)  d \sigma_{2 n - 1}  (\xi)}, \\ 
 \vt_{V} (a) & = & \lim_{r \to 0^{+}} \frac {\max_{\vert z - a  
 \vert = r} V (z)}{ \log r},
 \label{eq:fkiselman}
 \end {eqnarray}
 where $d \sigma_{2 n - 1}$ is the normalized area measure on the 
 unit sphere $\bd \B$.
 
 Observe that when $f$ is a holomorphic function near the point $a$ such that $f \not \equiv 0$ then $\vt_{\log \vert f\vert} (a)$ is the order of vanishing of $f$ at the point $a$. From the  formula (\ref {eq:fkiselman}), it follows immediately  
 that $\vt_{V} (a) = 0$ if $V (a) > - \infty.$ This formula shows 
 that the Lelong number $\vt_{V} (a)$ can be viewed as the weight 
 of the logarithmic singularity of $V$ at the point $a.$
 
 \section{The minimum principle for logarithmic potentials} 
 
 Recall that  the Lelong class on $\C^n$ is defined as follows 
 $$\mc L (\C^n) := \{v \in PSH (\C^n) ; v (z) \leq \log^+ \vert  
 z\vert + O (1), \forall z \in \C^n\}.$$
 To a function $V \in \mc L (\C^n)$ there is associated a Robin 
 function as follows (see [27], [37], [6]). For $z \in \C^n \sm \{0\},$ set
 $$\rho_V (\zeta) := \limsup_{\la \in \C, \la \to \infty} (V (\la 
 z) - \log \vert \la z\vert).$$
 Since this function is constant on any complex line of $\C^n$  
 passing through the origin, it follows that $\rho_V$ is a well 
 defined function on projective space $\mb P^{n - 1}$ which can be 
 viewed as the 
 hyperplane at infinity in $\C^n$. 
  Then following Bedford and Taylor ([6]), we introduce the 
 following class.
 $$\mc L_{\star} (\C^n) := \{ V \in \mc L (\C^n) ; \rho_V \not 
 \equiv - \infty \}.$$
 Let $\w_0$ be the Fubini-Study form on $\mb P^{n - 1}$ normalized 
 by the condition $\int_{\mb P^{n - 1}} {\w_0}^{n - 1} = 1.$ Then if $V \in 
 \mc L_{\star} (\C^n)$, the Robin function $\rho_V$ is an  
 $\w_0-$\psh on $\mb P^{n - 1}$,  in  the sense that it is 
 upper semi-continuous on $\mb P^{n - 1}$ and satisfies the 
 condition $dd^c \rho_V + \w_0 \geq 0$ in the sense of currents 
 on $\mb P^{n - 1}$. 
 
 An interesting fact concerning the class $\mc L_{\star} (\C^n)$ 
 is the following kind of Riesz representation formula which is  
 well known in one  variable but seems not to be known in $\C^n$.  
 This formula was obtained earlier with my student Fatima Amghad  
 but never published 
 (see [2]). We will give a proof below since we will use it later.
 \begin{lem} Any function $V \in \mc L_{\star} (\C^n)$ admits the 
 following representation formula
 \begin{equation} 
 V (z) = \int_{\C^n} \log \vert \zeta - z\vert dd^c V \we (dd^c 
 \log \vert \zeta - z \vert)^{n - 1} + \int_{\mb P^{n - 1}} \rho_V 
 \w_0^{n - 1},
 \label{eq:Lstar}
 \end{equation}
 for all  $ z \in \C^n.$
 \end{lem}
 \demo By translation we may assume that $z = 0$ is the origin in 
 $\C^n$. By the classical Poisson-Jensen formula, for $0 < r < R$, 
 we have
 $$ \int_{\vert \zeta \vert = 1} V (R \zeta) d \si_{2 n - 1} - 
 \int_{\vert \zeta \vert = 1} V (r \zeta) d \si_{2 n - 1} (\zeta)
 = \int_r^R \vt_V (t) \frac{d t}{t},$$
 where $\vt_V (t) := \vt_V (0,t)$ is the projective mass of the 
 current $dd^c V$ on the ball $\B (0,t)$.
 From this formula it follows that 
 $V (0) > - \infty$ iff $\int_0^R \vt_V (t) \frac{d t}{t} < + 
 \infty$, which implies that $\lim_{t \downarrow 0} \vt (t)  \log t  = 0.$  
 
 Assume first that $V (0) > - \infty$.  
  Then integration by part gives the following formula
 \begin{eqnarray*}
 \int_{\vert \zeta \vert = 1} V (R \zeta) d \si_{2 n - 1} - 
 \int_{\vert \zeta \vert = 1} V (r \zeta) d \si_{2 n - 1} (\zeta)
 & = & \vt_V (R) \log R - \vt_V (r) \log r \cr
 & & -  \int_r^R \log t \ d \vt_V (t),
 \end{eqnarray*}
 and then letting $r \downarrow 0$, we get the following formula
 \begin{equation}
 \int_{\vert \zeta \vert = 1} V (R \zeta) d \si_{2 n - 1} - V (0)
 = \vt_ V (R) \log R  -  \int_r^R \log t d \vt_V (t).
 \label{eq:FormRep}
 \end{equation}
 Now observe that by approximating $V$ by bounded functions $V_j := 
 \sup \{V,- j\}$ we see from the above formula that $V (0) = - 
 \infty$ iff $\int_0^R \log t \ d \vt_V (t) = - \infty$ and then  
 the formula (\ref{eq:Lstar}) holds also in this case.
 
 Therefore it is enough to prove the formula (\ref{eq:Lstar}) when  
 $V (0) > - \infty$.
 
 In that case, the formula (\ref{eq:FormRep}) yields the following 
 one
 $$\int_{\vert \zeta \vert = 1} V (R \zeta) d \si_{2 n - 1} - \vt_ 
 V (R) \log R = 
 V (0) -  \int_{\vert \zeta\vert < R} \log \vert\zeta\vert dd^c V 
 \we (dd^c \log \vert \zeta\vert)^{n - 1}.$$
 Now since $V \in \mc L_{\star} (\C^n)$, it follows that
 $$ \lim_{R \to + \infty} (\int_{\vert \zeta \vert = 1} V (R \zeta) 
 d \si_{2 n - 1} - \log R) = \int_{\vert \zeta \vert = 1} \rho_V 
 (\zeta)  d \si_{2 n - 1},$$
 and then  
 $$\lim_{R \to + \infty} (\log R)^{- 1} \int_{\vert \zeta \vert = 
 1} V (R \zeta) d \si_{2 n - 1} = 1.$$
 Therefore from Poisson-Jensen's formula it follows immediately that
 $\lim_{R \to + \infty} \vt_V (R) = 1,$  
 $\lim_{R \to + \infty} \int_{\vert \zeta\vert < R} \log 
 \vert\zeta\vert dd^c V \we (dd^c \log \vert \zeta\vert)^{n - 1}$ 
 is finite and then
 $\lim_{R \to + \infty} (1 - \vt_V (R)) \log R = 0$.
 
 Therefore we deduce the following formula
 $$\int_{\vert \zeta \vert = 1} \rho_V (\zeta)  d \si_{2 n - 1}
 = V (0) - \int_{\vert \C^n} \log \vert\zeta\vert dd^c V \we (dd^c 
 \log \vert \zeta\vert)^{n - 1}.$$
 On the other hand, using Fubini's theorem for the projection $\pi 
 :  \mb S_{2 n - 1} \lra \mb P^{n - 1}$, we deduce that
 $$ \int_{\vert \zeta\vert = 1} \rho_V (\zeta) d^c \log \vert 
 \zeta\vert \we (dd^c \log \vert \zeta\vert)^{n - 1} =
  \int_{\mb P^{n - 1}} \rho_V \w_0^{n - 1}.$$
  Now observe that the normalized area measure $\si_{2 n - 1}$ on 
 the unit sphere $\mb S_{2 n - 1}$ coincides with
 the restriction to $\mb S_{2 n - 1}$ of the form 
 $d^c \log \vert \zeta\vert \we (dd^c \log \vert \zeta\vert)^{n - 
 1}$, which implies the required formula. \fin
 
 It follows from this formula that the following class of 
 plurisubharmonic functions 
 $$\mc L_{og} (\C^n) := \{ V \in \mc L_{\star} (\C^n) ; \int_{\mb 
 P^{n - 1}} \rho_{V} \w_0^{n - 1} = 0\}$$
 is a natural generalization of the class of classical logarithmic 
 potentials. For this reason we will call this class the class of 
 logarithmic potentials in $\C^n$.
 
 Connected with this class there is a natural capacity on $\C^n$ 
 which is a generalization of the classical logarithmic capacity.
 For any subset $E \sub \C^n$, recall that the Siciak-Zahariuta's 
 extremal function associated to $E$ is defined by the following
 formula (see [32], [33], [36]).
 $$V_E (z) := \sup \{V (z) ; V \in \mc L (\C^n) ; \sup_E V \leq 
 0\}.$$
 Then the logarithmic capacity of the set $E$ will be defined by 
 the following formula
 $$ C_{\log} (E) :=  \exp (- \int_{\mb P^{n - 1}} \rho_{V_E^*} 
 \w_0^{n - 1}).$$
 Observe that if $E$ is pluripolar then $V_E^* \equiv + \infty$ and 
 then $C_{\log} (E) = 0$. On the other hand if $E \Sub \C^n$ is 
 non-pluripolar then $V_E^* \in \mc L (\C^n)$ and $\displaystyle 
 {\rho_{V_E^*}}$ is bounded on $\mb P^{n - 1}$, which proves that  
 $C_{\log} (E) > 0$.
 
 Moreover using a result concerning the convergence of Robin's functions from [6], it's possible to prove that 
 this logarithmic capacity is a Choquet capacity on $\C^n$. 
 
 Observe that this capacity is related to the class of 
 logarithmic potentials by the following formula
 \begin{equation}
  \log C_{\log} (E) = \inf \{ {\sup_{E}}^{\star} V ; V \in \mc L_{og} 
  (\C^n)\},
  \label{eq:caplog1}
  \end{equation}
 where ${\sup_{E}}^{\star} V := \inf\{ \sup_{E \sm A} V ;   A \subset E \  
 \mrm{is  \ pluripolar}\}$ is the quasi-essential upper bound of 
 $V$ on the set $E$.
 Indeed, it is easy to see from the definition of $V_E^{\star} $ and   
 results of [5] that $\sup_{E}^{\star}  V_E^{\star}  = 0$ and then
 $$ V_E^{\star} (z) = \sup \{V (z) ; V \in \mc L (\C^n) ; {\sup_{E}}^{\star} V = 
 0\}, \ \ z \in \C^n.$$
  Therefore we easily get the following formula
 \begin{eqnarray}
 - \log  C_{\log} (E) & = & \sup \{ \int_{\mb P^{n - 1}} \rho_V 
 \w_0^{n - 1} ; V \in \mc L (\C^n), {\sup_E}^{\star} V = 0\}
 \nonumber  \\
 & = & \sup \{ \int_{\mb P^{n - 1}} \rho_V \w_0^{n - 1}  - 
 {\sup_E}^{\star} V  ; V \in \mc L (\C^n)\},
 \label{eq:caplog2}
 \end{eqnarray}
 where the supremum is attained for $V = V_E^{\star} $.
 
 Now from the formula (\ref{eq:caplog2}), it follows that
 $$\log  C_{\log} (E) = \inf \{ {\sup_E}^{\star} V  - \int_{\mb P^{n -  
 1}} \rho_V \w_0^{n - 1} ; V \in \mc L (\C^n)\},$$
 which proves the formula (\ref{eq:caplog1}).
 
 Recall that  there is another constant which was classically  
 called the logarithmic capacity and defined  for a compact 
 subset $K \sub \C^n$ by the following formula ([33], [27], [22],  
 [37])
 $$ \tau (K) :=   \exp
 (- \limsup_{\vert z\vert \to + \infty} (V_K^{\star} (z) - \log \vert  
 z\vert) = \exp (- \sup_{\vert z\vert = 1} \rho_{V_K^{\star}}).$$
 From a well known inequality (see [1], [14], [33]), it follows  
 easily that there exists a constant $\kappa_n > 0$ such that
 $$ \tau (K) \leq C_{\log} (K) \leq \kappa_n \cdot \tau (K)$$
 for any compact subset $K \sub \C^n$.
 
 We want to prove the following generalization of the  
 Cartan-Boutroux lemma for the class $\mc L_{og} (\C^n)$ of  
 logarithmic potentials.
 \begin{thm} 
 For any real number $0 < \eta < 5$ and any function
 $V \in \mc L_{og} (\C^n)$, the following  lower bound  
 \begin{equation} 
 V (z) \geq  - \log (5 e \slash \eta),
 \label{eq:estimate-V}
 \end{equation}
 holds for all $z \in \C^n$, outside the union of a countable 
 family  of euclidean balls $(\B(z_j,r_j))$ of radii $(r_j)$ less 
 than $ \eta$, satisfying the following condition 
 \begin{equation}
 \sum_{\B (z_j,r_j) \cap \B_R \neq \emptyset} r_j^{2 n - 2 + \al} <  5^{2 n - 2} (R + \eta)^{2 n - 2} \eta^{\al} 
 \slash \al, \forall R > 0.
 \label{eq:estimateHC1}
 \end{equation}
 In particular, the exceptionnal set $E_{\eta} \sub \C^n$ where the lower  estimate (\ref{eq:estimate-V}) is not satisfied is a Borel set  
  for which the following estimate 
 \begin{equation}
 h^{2 n - 2 + \al}_{\eta} (E_{\eta} \cap \B_{R}) < 5^{2 n - 2} (R + \eta)^{2 n - 2} \eta^{\al} \slash \al, \ \forall R > 0,
 \label{eq:estimateHC2}
 \end{equation}
 holds.
 \end{thm}
 \demo Let $V \in \mc L_{og} (\C^n)$.
 First observe that considering the Stieltjes' integral with 
 respect to the increasing function $g: t \lto \vt_V (z,t)$, we can 
 write 
 \begin{equation}
 V (z) = \int_0^{+ \infty}  \log t \ d g (t) \geq   \int_0^{1}  \log t \ d g (t), \ \forall z \in \C^n.
        \label{eq:ineqV}
 \end{equation}
  Now fix the real numbers $0 < \al \leq 2$ and $0 < \vep < 1  
  \slash 5$ and let $A > 0$ be a real number to be specified 
 later in terms of $\vep$ and $\al$. Then denote by $G = 
 G_{\vep, \al}$ the subset of "good" points $z \in \C^n$ for which   we have the following bound on the projective mass
 $$\vt_V (z,t) \leq A t^{\al}, \ \forall 0 < t \leq \vep.$$
 Observe that this implies in particular that $\lim_{t \to 0} \vt_V (z,t) \log t = 0$ for  
 $z \in G$, which implies  that $V (z) > - \infty$ for $z \in G$  
 and thus the set $G$ doesn't meet the polar set of $V$.
 
 Now fix a point $z \in G$. Then $V (z) > - \infty$ and integration 
 by parts in the  Stieltjes' integral \label{eq:ineqV} implies  
 immediately that
 \begin{eqnarray}
 V (z) & \geq & - \int_0^1 \vt_V (z,t) \frac{d t}{t} \cr
   \nonumber
       & \geq & - \int_0^{\vep} \vt_V (z,t) \frac{d t}{t} -
     \int_{\vep}^1 \vt_V (z,t) \frac{d t}{t}  \cr
     \nonumber
  &\geq & - \frac{A}{\al} \vep^{\al} -  \int_{\vep}^1 \vt_V (z,t)  
  \frac{d t}{t}.
         \label{eq: F-estimate}
   \end{eqnarray}
   Since $V \in \mc L (\C^n)$ and due to our normalisation of the   
 complex Monge-Amp\`ere operator, we know that
   $$ \vt_V (z,t) \leq \int_{\C^n} dd^c V \we (dd^c  \log \vert z -  \zeta\vert)^{n - 1} \leq 1, \ \forall z \in \C^n, \forall t > 0.$$
 Therefore we derive the following estimate
  $$V (z) \geq - \frac{A}{\al} \vep^{\al} -  \log 
 (1 \slash \vep), \forall z \in  
 G.$$
 Choosing $A := \al \vep^{- \al}$, we obtain the 
 following lower bound
 $$ V (z) \geq  - \log (e \slash \vep), \ \forall z \in G.$$
 Therefore we obtain the required inequality if we can  estimate 
 properly the size of the exceptionnel set $E := \C^n \sm G$.
 
  From the definition of the set $E$, it follows that for any $z 
 \in E$ there exists a real number $0 < t_z < \vep$ such that 
 $$\vt_V (z,t_z) > A t_z^{\al}.$$
  On the other hand, an easy computation shows that  
 \begin{equation}
     r^{2 - 2 n} \mu_V (\B (z,r)) = \vt_V (z,r) \leq 1, \forall z  
 \in \C^n, \forall r > 0.
 \end{equation}
  Therefore  we get 
 $$\mu_V (\B (z,t_z)) \geq t_z^{2 n - 2} \vt_V (z,t_z) >  A 
 {t_z}^{2 n - 2 + \al}, \forall z \in E.$$
  We want to exclude such balls. Since $(\B (z,t_z))_{z \in E}$ is  
 a covering of the exceptionnal set $E$ by open euclidean balls,
 by a Vitalli type $5-$covering lemma (see [29]), there exists a 
 countable subfamily of mutually disjoint balls $(\B (z_j,t_j))_{j 
 \in \N}$ such 
 that the corresponding $5-$family of balls $(\B (z_j,5 t_j))_{j  
 \in \N}$  
 covers $E$. 
  
 Now fix $R > 0$ and consider the familly $(\B (z_j,5 t_j))_{j  
 \in J_R}$ of those balls which  
 intersect  the ball $\B_R$. Then we obtain the following estimate
 $$ A \sum_{j \in J_R} {(5 t_j)}^{2 n - 2 + \al} < 5^{2 n - 2 + \al} \sum_{j \in J_R} 
 t_j^{2 n - 2} \vt_V (z_j, t_j)) \leq  5^{2 n - 2 + \al}  
 \sum_{j \in J_R} \mu_V (\B (z_j,t_j)).$$
 Observing that 
 $$ \sum_{j \in J_R} \mu_V (\B (z_j,t_j)) =  \mu_V (\bigcup_{j \in J_R} \B (z_j,t_j) 
 \leq  \mu_V (\B (0,R + \vep)) \leq (R + \vep)^{2 n - 2},$$
  we conclude that
 $$ \sum_{j \in J_R} {(5 t_j)}^{2 n - 2 + \al} < (5 R + 1)^{2 n -  
 2} 5^{\al} \slash A = (5 R + 1)^{2 n - 2} (5 \vep)^{\al} \slash 
 \al.$$
 Taking $r_j := 5 t_j$ and $ \vep := \eta \slash 5$, we obtain the 
 theorem  since the family of balls $(\B (z_j,r_j))_{j \in J_R}$ 
 covers $E \cap \B_R$.\fin
 Observe that the above result yields immediately a precise 
 estimate on the ${2n - 2 + \al}-$ Hausdorff contents of the 
 plurisubharmonic lemniscates associated to functions in the class 
 $\mc L_{\star} (\C^n)$.
 \begin{cor} Let $V \in \mc L_{\star} (\C^n)$ and $0 < \vep < 1 
 \slash e$. Then for any $\alpha \in ]0, 2]$, the ${2n - 2 + 
 \al}-$Hausdorff content of the associated plurisubharmonic 
 lemniscate $ E (V,\vep) := \{ z \in \C^n ; V (z) \leq \log \vep 
 \}$ satifies the following estimates 
 $$ h^{2 n - 2 + \al}_{\eta} (E_{\vep} \cap \B_{R}) < 5^{2 n - 2} 
 (R + 5)^{2 n - 2} (5 e \vep)^{\al} \slash \al, \ \forall R > 0,$$
 holds.
 \end{cor}
 
 From the theorem it's also possible to deduce the following  
 comparison inequality between certain relative Hausdorff contents  
 and the logarithmic capacity defined above.
 
 \begin{cor} For any real number $0 < \al \leq 2$  and any subset  
 $K  \sub \B$, we have
 $$ h^{2 n - 2 + \al} (K) \leq \frac{c_n}{\al} \bigl(5 e C_{\log}  
 (K)\bigr)^{\al},$$
 where $c_n = 5^{2n - 2} (1 + 1 \slash e)^{2n - 2}.$
 \end{cor}
 \demo  Recall from the formula (\ref{eq:caplog1}) that 
 $$\log C_{\log} (K) = \inf \{ \sup_K^{\star} V \, ; \, V \in \mc  
 L^+ (C^n) ; \int_{\mb P^{n - 1}} \rho_V \, {\w_0}^{n - 1} = 0 \},$$
 where $\sup_K^{\star}  V$ is the quasi-essential upper bound of  
 $V$ on $K$.
  
 Assume first that $ K \sub \B_r$ with $r := 1 \slash e$ so that $  
 C_{\log} (K) < 1 \slash e.$ Then let  $c$ be an arbitrary real  
 number such that $ C_{\log} (K) < c < 1  \slash e$. Then there  
 exists a function $V \in \mc L_{og} (\C^n)$ and a pluripolar  
 subset $A \sub E$ such that $\sup_{K \sm A} V < c $ so that $K \sm  A \sub K_{c} := \{ z \in \B_r ; V (z) <  \log c\}$. 
  By the minimum principle for the class $\mc L_{og} (\C^n)$, with  
 $\eta = 5 e c <  5$, we have $V  (z) \geq - \log (5 e \slash \eta)  = \log c$ for $z \in \B \sm E$,   
 where $E \sub \B$ is a Borel set satisfying $h^{2 n - 2 + \al}   
 (E) \leq 5^{2 n - 2} (r + \eta \slash 5)^{2 n - 2} \eta^{\al}  
 \slash \al = 5^{2 n - 2} (r + 1)^{2 n - 2} (5 e c)^{\al} \slash  
 \al$. Since by definition $K \sm A \sub K_{c} \sub E$, we deduce  
 that 
 $h^{2 n - 2 + \al} (K \sm A) \leq (5 r + 1)^{2 n - 2} (5 e c)^{\al} \slash \al.$
 
 Since $c > C_{\log} (K)$ is arbitrary, we obtain the 
 inequality 
 $$h^{2 n - 2 + \al} (K \sm A) \leq (5 r + 1)^{2 n - 2} (5 e C_{\log} (K))^{\al} \slash \al.$$
 Since $A \in \C^n$ is pluripolar, it is polar in $\R^{2 n}$ and then from a well known result in classical potential theory  (see [23]), it follows that $h^{2 n - 2 + \al} (A) = 0$ and then $h^{2 n - 2 + \al} (K) = h^{2 n - 2 + \al} (K \sm A),$ which proves the required estimate in the case when $K \sub B_r$.
  
 Now if $K \sub \B$ is any subset, it is enough to apply the last inequality to the  set $r \cdot K \sub r \B$ to obtain the required inequality.
 \fin
 \noindent {\bf Remarks:} 1)  Observe that if $n = 1$ and $\al = 
 1$, our estimate \ref{eq:estimateHC2}) reduces precisely to that given in [26]. This extends the Cartan-Boutroux lemma  
 except that the constant in our estimate is $5 e$ instead of $2 
 e$. \\ 
 2) When $n \geq 2$ our estimate shows how the relative  
 Hausdorff content of the exceptionnal set with respect to the ball  $\B_R$ is asymptotically small when $R \to + \infty$.\\
 3) Observe that for $\al = 2$ similar estimates  in terms of  
 the relative logarithmic capacity was obtained in [8].
 
 \section{A lower bound for plurisubharmonic functions} 

 Let $\B$ be the euclidean unit open ball in $\C^n.$ For each $z  
 \in \B,$ we denote by $\Phi_{z}$ the involutive automorphism of  
 the unit ball $\B$ which takes the point $z \in \B$ to the origin. 
 Then the pluricomplex Green function $G_z (\zeta) := G (\zeta,z)$  
 of the unit ball $\B$ with a logarithmic pole at the point $z \in   \B$ is given by the formula
 \begin{equation}
  G_{z} (\zeta) := \log \vert \Phi_{z} (\zeta)\vert  \ ,
  \ (z,\zeta) \in \B \times \B.
        \label{eq:green}
 \end{equation}
 It is easy to see that the following fundamental Monge-Amp\`ere   
 equation 
 \begin{equation}
  (dd^c G_{z})^n = \delta_{z}
  \label{eq:f-equation}
 \end{equation}
 holds in the sense of currents on $\B,$ where $\delta_{z}$ 
 is the unit Dirac mass at the point $z.$  \\ 
 It is well known that the formula $ d_{\B} (z,\zeta) := \vert   
 \Phi_z (\zeta)\vert$ defines a distance on the unit ball $\B$  
 which is related to the Bergman distance $\rho_{\B}$ by the following 
 formula
 $$d_{\B} (z,\zeta) = \tanh \frac{\rho_{\B} (z,\zeta)}{\sqrt{n + 1}}.$$
 Now  consider  the corresponding ball of center $z \in \B$ and 
 radius $r  \in ]0 , 1[$ defined by
 $$\omega_{z} (r) := \{\zeta  \in \B \ ; \ \vert 
 \Phi_{z} (\zeta) \vert <  r \},$$
 and define the following "invariant projective mass" function
 \begin{equation}
  \theta_{V} (z , r) := \int_{\omega_{z} (r)} 
  dd^c V \wedge (dd^c G_{z})^{n - 1},  
        \label{eq:projm}
 \end{equation}
 for $z \in \B$ and $0 < r < 1$ and observe that 
 \begin{equation}
  \theta_{V} (z , r) = \int_{\B_r} 
 dd^c V \circ \Phi_z \wedge (dd^c \log \vert \zeta \vert)^{n - 1}   
 = \vt_{V \circ \Phi_z} (0,r),
        \label{eq:projequ}
 \end{equation}
 and then 
 $\lim_{r\to 0} \theta_{V} (z , r) = \vt_{V \circ \Phi_z} (0) = 
 \vt_V (z)$ for any $z \in \B$ since $\Fi_z$ is an automorphism 
 taking the origin to the  point $z$ (see [15], [19]).  

 First we prove the following lemma which is similar to a result  
 of H. Milloux concerning a lower estimate for monic polynomials 
 with all zeros in the unit disc, the exceptional set being  
 estimated in terms of non euclidean distance on the unit disc (see  [30]).
 \begin{lem} Let $V$ be a \psh function on  the euclidean open unit  ball $\B \sub \C^n$ with bounded Riesz mass $\mu_V (\B) := (1  
 \slash 2 \pi) \int_{\B} \Delta V < + \infty$. Let us define the 
 following pluricomplex Green potential 
 \begin{equation} 
 \mc G_{V} (z) := \int_{\B} G_z dd^c V \we (dd^c G_z)^{n - 1}, \ z  
 \in \B. 
 \label{eq:Gpotential}
 \end{equation} 
 Then there exists a constant $c_n > 0$ such that for any real  
 numbers $0 < s < 1$  and $0 < \eta < \min \{3 s,1\},$ the following  lower bound  
 \begin{equation} 
 \mc G_{V} (z) \geq  - \theta_V (z,s) \log (3 \slash \eta) - c_n \  
 \mu_V (\B) \log (e \slash s),
 \label{eq:estimate-GV}
 \end{equation}
 holds for all $z \in  \B$, outside the union of a countable family  of pseudo-balls $(\w_{z_j} (r_j))$ of radii $(r_j)$ less than $  
 \eta$ and satisfying the condition 
 $\sum_j r_j^{2 n - 2 + \al} < 9^{n - 1} \eta^{\al} \slash \al$. 

 In particular, the exceptionnal set $E \sub \B$ where the lower  
 bound (\ref{eq:estimate-GV}) does not hold is a Borel set whose  
 invariant $\eta-$Hausdorff content of dimension $2 n - 2 + \al$ 
 satisfies the estimate $\ti{h}^{2 n - 2 + \al}_{\eta} (E) < 9^{n -  1} \eta^{\al} \slash \al$.
 \end{lem}
 \demo
 To estimate the function $\mc G_{V}$ given by the formula   
 (\ref{eq:Gpotential}), observe that by considering the Stieltjes' integral  with respect to the increasing function $g: t \lto \theta_V (z,t)$,  
 we can write the formula (\ref{eq:Gpotential}) as follows
 \begin{equation}
 \mc G_{V} (z) =  \int_0^1  \log t \ d g (t), \ \forall z \in \B.
        \label{eq:tiV}
 \end{equation}
 Now fix the real numbers $0 < \al \leq 2$ and $0 < \vep < \min  
 \{s, 1 \slash 3\}$ and let $A > 0$ be a real number to be  
 specified later in terms of $\vep$ and $\al$. Then denote by $U = 
 U_{\vep, \al}$ the  subset of "good" points $z \in \B$ for which  
 we have the following bound on the invariant projective mass
 $$\theta_V (z,t) \leq A t^{\al}, \ \forall 0 < t \leq \vep.$$
 Observe that this implies in particular that $\vt_V (z) = 0$ for  
 $z \in U$, which shows that the set $ U$ doesn't  
 contain the logarithmic singularities of $V$.
 
 Now fix a point $z \in U$. Then integration by parts in the  
 Stieltjes' integral \label{eq:tiV} implies immediately that
 \begin{eqnarray}
 \mc G_{V} (z) & = & - \int_0^1 \theta_V (z,t) \frac{d t}{t} \cr
   \nonumber
       & = & - \int_0^{\vep} \theta_V (z,t) \frac{d t}{t} -
     \int_{\vep}^1 \theta_V (z,t) \frac{d t}{t}  \cr
     \nonumber
  &\geq & - \frac{A}{\al} \vep^{\al} -  \int_{\vep}^1 \theta_V (z,t)  
  \frac{d t}{t}.
         \label{eq: F-estimate}
   \end{eqnarray}
 On the other hand, we can write
 $$ \int_{\vep}^1 \theta_V (z,t) \frac{d t}{t} =  \int_{\vep}^{s}  
 \theta_V (z,t) \frac{d t}{t} +  \int_{s}^1 \theta_V (z,t) \frac{d  
 t}{t}.$$
 Now observe that 
 $$ \int_{\vep}^s \theta_V (z,t) \frac{d t}{t} \leq \theta_V (z,s) \log  
 (s \slash \vep),$$
 and
 $$ \int_{s}^1 \theta_V (z,t) \frac{d t}{t} \leq  \theta_V (z,1) \log (1  
 \slash s).$$
 Therefore we conclude that
 $$\mc G_{V} (z) \geq - \frac{A}{\al} \vep^{\al} - \theta_V (z,s) \log 
 (s \slash \vep) - \theta_V (z,1) \log (1 \slash s), \forall z \in  
 U.$$
 Choosing $A := \al c_n \ \mu_V (\B) \vep^{- \al}$, we obtain the 
 following lower bound
 $$ \mc G_{V} (z) \geq - c_n \ \mu_V (\B) - \theta_V (z,s) \log (s 
 \slash \vep) -
 c_n \ \mu_V (\B) \log (1 \slash s), \ \forall z \in U.$$
 
  Therefore we obtain the required inequality if we can  estimate 
  properly the size of the exceptionnel set $E := \B \sm U$.
 
   From the definition of the set $E$, it follows that for any $z 
 \in E$ there exists a real number $0 < t_z < \vep$ such that 
 $$\theta_V (z,t_z) > A t_z^{\al}.$$
  On the other hand, an easy computation shows that there exists a 
 constant $c_n > 0$ such that 
 \begin{equation}
 \theta_V (z,r) \leq c_n \  r^{2 - 2 n} \mu_V (\w_z (r))),
  \forall z  \in \B, \forall r \in ]0, 1 \slash 3[.
 \end{equation}
  Therefore we get 
 $$c_n \mu_V (\w_z (t_z)) \geq t_z^{2 n - 2} \theta_V (z,r_z) >  A 
 {t_z}^{2 n - 2 + \al}, \forall z \in E.$$
  Since $(\w_z (t_z))_{z \in E}$ is a covering of the exceptionnal 
 set $E$ by invariant pseudo-balls, by an elementary Vitalli type 
 $3-$covering  
 lemma for the ball with the hyperbolic distance (see [31], [34]), 
 there exists a countable subfamily of mutually disjoint 
 pseudo-balls $(\w_{z_j} (t_j))$ such that the corresponding 
 $3-$family of pseudo-balls $(\w_{z_j} (3 t_j))$ covers $E$. 
  Moreover we obtain the following estimate
 $$ A \sum_j {(3 t_j)}^{2 n - 2 + \al} < 3^{2 n - 2 + \al} \sum_j  
 t_j^{2 n - 2} \theta_V (z_j, t_j)) \leq  3^{2 n - 2 + \al} c_n  
 \sum_j \mu_V (\w_{z_j} (t_j)).$$
 Observing that  
 $ \sum_j \mu_V (\w_{z_j} (t_j)) = \mu_V (\bigcup \w_{z_j} (t_j)) 
 \leq  \mu_V (\B)$, we conclude that
 $$ \sum_j {(3 t_j)}^{2 n - 2 + \al} < 9^{n - 1} 3 ^{\al} c_n  
 \mu_V (\B) \slash A = 9^{n - 1} (3 \vep)^{\al} \slash \al.$$
 Taking $r_j := 3 t_j$ and $ \vep := \eta \slash 3$, we obtain the 
 theorem  since the family of pseudo-balls $(\w_{z_j} (r_j))$ 
 covers $E$.\fin
 
 As a consequence of the last lemma we obtain a minimum principle for the Cegrell's class $\mc F (\B)$ ([12]) on the unit ball in terms of the pseudo-distance on the unit ball.
 \begin{prop} Let 
 Then there exists a constant $c_n > 0$ such that for any real  
 numbers $0 < \eta < 1$ and any function $\vfi \in \mc F (\B)$ with $\int_{\B} (dd^c \vfi)^n \leq 1$, the following lower bound  
 \begin{equation} 
 \vfi (z) \geq  -  \log (c_n \slash \eta)  
 \label{eq:estimate-GV}
 \end{equation}
 holds for all $z \in  \B$, outside the union of a countable family  of pseudo-balls $(\w_{z_j} (r_j))$ of radii $(r_j)$ less than $  
 \eta$ and satisfying the condition 
 $$\sum_j r_j^{2 n - 2 + \al} < 9^{n - 1} \eta^{\al} \slash \al.$$ 

 In particular, the exceptionnal set $E \sub \B$ where the lower  
 bound (\ref{eq:estimate-GV}) does not hold is a Borel set whose  
 invariant $\eta-$Hausdorff content of dimension $2 n - 2 + \al$ 
 satisfies the estimate 
 $$\ti{h}^{2 n - 2 + \al}_{\eta} (E) < 9^{n -  1} \eta^{\al} \slash \al.$$
 \end{prop} 
 \demo From the definition of the class $\mc F (\B)$ (see [12]), it  follows that there exists a decreasing sequence of \psh functions   on $\B$ with boundary values $0$ which converges to $\vfi$ and  
 satisfies the condition $\sup_{j} \int_{\B} (dd^c \vfi_j)^n < + \infty$. 
 Applying the Poisson-Jensen formula to every $\vfi_j$ and taking  
 the limit we obtain the formula 
 $\vfi (z) = \mc G_{\vfi} (z)$ for $z \in \B$. 
 
 It's well known that functions from the class $\mc F (\B)$ have bounded Riesz mass  (see [12]) and then we can apply the last lemma to conclude that the estimate (\ref{eq:estimate-GV}) holds for $V = \vfi$ and $s = 1$.
 Now oberve that $\mu_{\vfi} (\B) \leq 1$ and $\theta_{\vfi} (z,1) \leq 1$ if  $\vfi \in \mc F (\B)$ and $\int_{\B} (dd^c \vfi)^n \leq 1$ (see [Ce]), which implies the required estimates. \fin
 
 Now let us state the main result of this section. For a fixed $0 < \rho < 1$,  
 define the following constant
 $$ \ka_n (\rho) := \frac{(1 + \rho)^n}{(1 - \rho)^n}.$$
 \begin{thm} Let $0 < \rho < 1$ be a real number. 
 Then there exists a positive constant $c_n (\rho) > 1$ such that  
 for any real numbers $0 < s < 1$, $0 < \eta < \min \{s , 1  
 \slash 3\}$ and any function $V$ \psh on a neighbourhood of the  
 closed euclidean unit ball $\ove{\B} \sub \C^n$ such that $V\leq  
 0$ on $\B$, the following lower bound 
 \begin{equation} 
 V (z) \geq  \kappa_n (\rho) \int_{\bd \B} V d \si_{2 n - 1} -  
 \theta_V (z,s) \log (3 \slash \eta) - c_n (\rho) \ \mu_V (\B) \log (e  \slash s), 
 \label{eq:estimate-u}
 \end{equation}
 holds for all $z \in \B_{\rho}$, outside the union of a countable  
 family of pseudo-balls $(\w_{z_j} (r_j))$ of radii $(r_j)$ not  
 exceding $\eta$ and satisfying the condition 
 $\sum_j r_j^{2 n - 2 + \al} < 9^{n - 1} \eta^{\al} \slash \al$. 

 In particular, the exceptionnal set $E \sub \B_{\rho}$ where the  
 lower bound (\ref{eq:estimate-u}) does not hold is a Borel set 
 whose  Hausdorff content satisfies the estimate $h^{2 n - 2 + 
 \al}_{\eta} (E) < 9^{n - 1} \eta^{\al} \slash \al$.
 \end{thm}
 \demo 
 By the Jensen-Poisson-Szeg\"o formula, we get the following  
 representation formula 
 \begin{equation}
  V (z)  =  \int_{\partial \B} V  \ d^c G_{z} \wedge 
 (dd^c G_{z})^{n - 1} +  \int_{\B} G_{z} \ dd^c V \wedge 
        (dd^c G_{z})^{n - 1}, 
  \label{eq:JL}
 \end{equation}
 for $\ z \in \B$ (see [31], [34]).  
 Recall that this formula follows easily from the fundamental 
 equation (\ref{eq:f-equation}) and the fact that 
 $G_{z}$ has boundary values $0.$ 
 
  Now let us write $V = \mc P_V + \mc G_{V}$ on $\B$, where 
 \begin{equation}
 \mc P_V (z) :=    \int_{\partial \B} V   d^c G_{z} \wedge 
        (dd^c G_{z})^{n - 1}, \ z \in \B,
        \label{eq:PV}
 \end{equation}
 and 
 \begin{equation}
 \mc G_{V} (z) := \int_{\B} G_{z} \ dd^c V \wedge 
   (dd^c G_{z})^{n - 1}, \ z \in \B.
   \label{eq:GV}
 \end{equation}
 It is well known that 
 $$d^c G_{z} \wedge (dd^c G_{z})^{n - 1} = \mc P (z,.) \ d \si_{2 n  - 1},$$
 where 
 $$ \mc P (z,\zeta) := \frac{(1 - \vert z\vert^2)^n}{\vert 1 -  
 z\cdot \ove{\zeta}\vert^{2 n}}, (z,\zeta) \in \B \times \bd \B,$$
 is the Poisson-Szeg\"o kernel of the open unit ball $\B$ and $d \si_{2 n - 1}$ is the normalized area measure on the unit sphere $\mb S_{2 n - 1} = \bd \B$.

 Therefore, since $V \leq 0$ on $\B$, it follows that the function 
 $\mc P_{V}$ satisfies the following lower  estimate
 \begin{equation}
  \mc P_{V} (z) \geq \ka_n (\rho)   \int_{\partial \B} V (\zeta) d  
  \si_{2 n - 1} (\zeta), \ \ \mrm{for} \ \ 
        \vert z \vert \leq \rho.
        \label{eq:minorGV}
 \end{equation}
 Then the estimate of the theorem follows from the lemma.\fin 
 \section{A minimum principle for plurisubharmonic functions}
  Here we will give a general version of the minimum principle for compact classes of \psh functions and derive a kind of {\it three-circle  minimum principle} for arbitrary plurisubharmonic functions. 
 \begin{thm}  Let $\W \sub \C^n$ be an open set, $K \sub \W$ a  
 compact set and let $\mc U \sub  PSH (\W)$ be a compact class.  Then  $\vt := \sup \{\vt_{u}(z) ; u \in \mc U, z \in K\} < + \infty\}$ and for any $\nu > \vt$ there exists   
 constants $C = C (K,\nu) > 0$ and $0 < \eta_0 < 1$ such that for any real number $0 < \eta \leq \eta_0$, any real number $0 < \al \leq 2$ and any function $u \in \mc U$ the following lower 
 bound
 $$u (z) \geq - \nu \log (C \slash \eta),$$
 holds for any $z \in K$ outside the union of a countable family of 
 balls of radii $(r_j)$ less than $\eta$ which satisfy the following condition  
 $$\sum_j r_j^{2 n - 2 + \al} < N 9^{n - 1} \eta^{\al} \slash \al,$$ 
 where $N$ is an integer depending only on $(K,\W)$.
 \end{thm}
 \demo By upper semi-continuity of Lelong numbers, there exits a real number $0 < s < 1$ small enough and an open neighbourhood $\w \Sub \W$ of $K$ such that $\theta_{u}(z,s) 
 < \nu$ for any $z \in \w$ and any $u \in \mc U$ (see [39]). Now take a finite number of concentric euclidean balls $B'_i \Sub \B_i \Sub \w$ ($1 \leq i \leq N$) such that the balls $(\B'_i)_{1 \leq i \leq N}$ cover $K$. By Theorem 5.2, for each $i$ there exits constants $\ka_i > 0$ and $c_i > 0$ such that the following lower bound 
 $$ V (z) \geq  \kappa_i \int_{\bd \B_i} V d \si_{2 n - 1} - \theta_V (z,s) \log (3 \slash \eta) - c_i \ \mu_V (\B_i) \log (e \slash s), $$
 holds for all $z \in \B'_i \sm E_i$, where $E_i \sub \B'i$ is a Borel set with $h^{2 n - 2 + \al}_{\eta} (E_i) < 9^{n - 1} \eta^{\al} \slash \al$.
 By compactness of $\mc U$ it follows that $ \int_{\bd \B_i} V d \si_{2 n - 1}$ and $ \mu_V (\B_i)$ are uniformly bounded for $V \in \mc U$. 
 Therefore taking $E = \bigcup_{1 \leq i \leq N} E_i$, we obtain for a suitable constant $C > 0$ the following lower  bound
 $$V (z) \geq - \theta_V (z,s) \log (3 \slash \eta) - C \geq \nu \log (3 \slash \eta) - C, $$
  for any $z \in K \sm E$.
 Since $h^{2 n - 2 + \al}_{\eta} (E) \leq \sum_i h^{2 n - 2 + \al}_{\eta} (E_i) < N 9^{n - 1} \eta^{\al} \slash \al$, we obtain the estimate of the theorem with an appropriate constant. \fin
 
 Recall that   the class $\mc U_\B := \{ v \in \mc L (\C^n) ; \sup_{\ove \B} v = 0\}$ is a compact class of  \psh functions. Applying the same proof as in the last theorem we obtain the following result which will be called the minimum principle for the Lelong class. 
 \begin{thm} For any $\nu > 1$ there exists   
 constants $C > 0$ and $0 < \eta_0 < 1$ such that for any real number $0 < \eta \leq \eta_0$, any real number $0 < \al \leq 2$, any real number $R \geq 1$ and any $V \in \mc L (\C^n)$, the following lower bound
 $$V (z) \geq \max_{B_R} V - \nu \log (C \slash \eta),$$
 holds for any $z \in \B_R$ outside the union of a countable family  of euclidean balls of radii $(r_j)$ less than $\eta R$ which satisfy the following condition  
 $$\sum_j r_j^{2 n - 2 + \al} <  9^{n - 1} R^{2 n - 2 + \al} \eta^{\al} \slash \al.$$ 
 \end{thm}
 \demo It is enough to prove the estimate for $R = 1$. Then applying the same method as in the proof of the previous theorem to the compact class $\mc U_{\B}$ and observing that for any $v \in \mc L (\C^n)$, $u := v - \sup_{\ove \B} v \in \mc U_{\B}$ we obtain the required lower bounds since in this cas $N = 1$. \fin  
 
 We want to apply the minimum principle for compact classes to  
 prove the following general result which extends the classical 
 minimum principle stated in section 2. 
 
 To motivate our result, let us recall that the maximum principle for plurisubharmonic functions can be used to prove the classical Hadamard three-circle inequality which can be stated as follows. Let $V$ be a \psh function on the euclidean ball $\B_R$ of radius $R > 0$ and let 
 $\si, \tau \in ]0 , 1[$ be real numbers with $0 < \si \leq \tau < 1$. Then the following upper bound
 $$  V (z) \ \leq  \ \sup_{\B_{\si R}} V + \rho (\si,\tau) \ \Bigl(\sup_{\B_{R}} V - \sup_{\B_{\si R}} V\Bigr),$$
 holds for any $z \in \B_{\tau R}$, where 
 $$\rho (\si,\tau) := \frac{\log (\tau \slash \si)}{\log (1 \slash \si)}.$$
 This can be viewed as a "three-circle maximum principle".
 
 Here we want to establish a kind of "three-circle minimum  
 principle" which is dual in some sense to the previous one.
 
 Define for $\si, \tau \in ]0 , 1[$ the following constant
 $$ \nu (\si,\tau) := \frac{1}{ \log \Bigl(\frac{1 + \si \tau}{\si  
 + \tau}\Bigr)}.$$
 \begin{thm} Let $\si, \tau \in ]0 , 1[$ be given real numbers.
 Then for any real number $\nu > \nu (\si,\tau)$ there exists   
 constants $C = C (\si,\tau,\nu) > 0$ and $0 < \eta_0 < 1$ such 
 that for any $R > 0$, for any \psh function
 $V$ on the euclidean ball $\B_R$ and for any $\al \in ]0 , 2]$ and 
 any $\eta \in ]0,\eta_0[$,
 the following lower bound 
 \begin{equation}
  V (z) \ \geq \ \sup_{\B_{\si R}} V + \nu \ \log (\eta \slash C) \ 
  \Bigl(\sup_{\B_R} V - \sup_{\B_{\si R}} V\Bigr),
 \label{eq:lb}
 \end{equation}
 holds for all $z \in \B_{\tau R}$ outside the union of a countable 
 family of euclidean balls of radii $(r_j)$ not exeeding $\eta R$ 
 and satisfying the estimate
 \begin{equation}
 \sum_j r_j^{2 n - 2 + \al} < 9^{n - 1} (R e)^{2 n - 2 + \al}  
 \eta^{\al} \slash \al.
 \label{eq:hc1}
 \end{equation}
 \end{thm}
 Observe that the theorem gives a precise uniform bound of the 
 relative Hausdorff content of dimension $(2n - 2 + \al)$ of the 
 exceptionnal set where the lower estimate (\ref{eq:lb}) does not 
 hold.  
 \demo Since the inequality (\ref{eq:lb}) and the condition  
 (\ref{eq:hc1}) are invariant under any homothetic map, we can 
 assume that  $ R = 1.$ Denote by $\mc U$ the class of \psh 
 functions on $\B$ such that $u  \leq 1$ and $ \max_{\B_{\si}} u  
 \geq 0$. Then $\mc U$ is a compact class of  
 \psh functions. To apply Theorem 6.1, we need to estimate the 
  Lelong numbers of the class $\mc U$ on the compact ball 
 $\ove{\B}_{\tau}$. 
  Indeed, if $u \in \mc U$  and $z \in \ove{\B}_{\tau}$, then if $0  < r < 1$, we have 
 $$\vt_{u} (z) = \vt_{u \circ \Phi_z} (0) \leq \frac{\sup_{\B} u 
 \circ \Phi_z - \sup_{\B_r} u \circ \Phi_z}{\log (1 \slash r)} \leq
 \frac{ 1 - \max_{\Phi_z (\B_r)} u }{\log (1 \slash r)}.$$
 Now it follows from a simple computation (see [30]) that for $z \in \C^n$ with $\vert z \vert \leq \tau$, the set $\Phi_z (\B_r)$ contains the ball $\B_{\si}$ precisely when  $r = (\si + \tau)\slash (1 + \si \tau),$ which implies that $\max_{\Phi_z (\B_r)} u \geq 0$.
 
 Then for any $u \in \mc U$ and any $z \in \ove{\B}_{\tau}$ we 
 have
 $$\vt_{u} (z) \leq \nu (\si,\tau).$$ 
 This inequality was first obtained in [7].
 
 Therefore from Theorem 6.1, it follows that for any $\nu > \nu (\si,\tau)$, there exists a 
 constant  $C > 0$ and a real number $\eta_0 \in]0,1[$ small 
 enough, for any $0 < \eta < \eta_0$, for any $\al \in ]0,2]$ and  
 for any $u \in \mc U$, we get the following uniform lower bound
 \begin{equation}
 u (z) \geq - \log (C \slash \eta),
 \label{eq:lbu}
 \end{equation} 
 for any $z \in \B_{\tau}$ outside the union of a countable family 
 of balls with radius $(r_j)$ less that $\eta$ and satistying the 
 estimate
 $\sum_j r_j^{2n - 2 + \al} < 9^{n - 1} \eta^{\al} \slash \al.$
 
 Now let $V$ be an arbitrary non constant \psh function on $\B$ so that $\sup_{\B} V > \sup_{\B_{\si}} V$.  
 Then the following function 
 $$u := \frac{V - \sup_{\B_\si} V}{\sup_{\B} V - \sup_{\B_\si} V}$$
 belongs to the class $\mc U$. Therefore this function satisfies 
 the lower bound (\ref{eq:lbu}) which yields the required lower 
 bound for $V$.  \fin
 
 As a consequence we give the following generalisation of the one  
 variable minimum principle stated in section 2.
 \begin{cor} There exists a constant $C > 1$ and a real number $0 <  \eta_0 < 1$ such that  for any real number $0 < \eta \leq \eta_0$ 	 and any \psh function  $V$ is on any euclidean  
 ball $\B_{\ti R} \sub \C^n$, where $\ti R := 2 e R$ and $R > 0$  
 satisfying the condition $V (0) = 0$,  the following lower bound
 \begin{equation}
  V (z) \geq -  \log (C \slash \eta) \max_{\B_{\ti R}} V  
  \label{eq:mp2}
  \end{equation}
 holds for any $z \in \B_{R}$ outside the union of a countable 
 family of balls of radii $(r_j)$ less than $\eta R$ which satisfy  
 the following condition
 \begin{equation}
 \sum_j r_j^{2 n - 2 + \al} < 9^{n - 1} (R e)^{2 n - 2 + \al}  
 \eta^{\al} \slash \al.
 \label{eq:hc2}
 \end{equation}
 \end{cor}
 \demo It is enough to apply the last result with $\ti R$ instead 
 of $R$, $\tau = 1 \slash (2 e)$ and $\si > 0$ small enough so that 
 $\nu (\si,\tau) < 1$, which is possible since $\nu (0,\tau) = 1 \ 
 \slash (\log (2 e)) < 1$. \fin
 Let us give an asymptotic formulation of the last result which  
 explains how the above three-circle minimum principle is dual to 
 the  three-circle maximum principle .
 
 Let us define for $\vep > 0$ the $(h_p,\vep)-$essentiel lower   
 bound of a \psh function $u$ on the ball $\B_{\tau R}$ as follows
 $$\mc I_{p,\vep} (u,\B_{\tau R}) := \sup \{ \inf_{\B_{\tau R} \sm  
 E} u ; E \sub \B_{\tau R}, {h^p}^* (E) \leq \vep\},$$
 where ${h^p}^*$ is the exterior Hausdorff content or order $p$.
 Since $u$ may have poles, the real numbers $\mc I_{p,\vep} (u,\B_{\tau R})$ may decrease to $- \infty$ when $\vep \downarrow  
 0$. 
 The minimum principle above gives the rate of convergence to $-  
 \infty$ of this minimum uniformly when the function $u$ is  
 suitably normalized. More precisely, we deduce the following  
 result.
 \begin{cor} Let $\mc U$ the class of \psh function on tha ball  
 $\B_R$ such that $u \leq 1$ on $\B_R$ and  $\max_{\si B} u = 0$.  
 Then for any $\al \in]0, 2]$, the following uniform asymptotic  
 estimate 
 $$\limsup_{\vep \downarrow 0} \Bigl(\sup_{u \in \mc U} \frac{\mc  
 I_{p,\vep} (u,\B_{\tau R})}{\log \vep}\Bigr) \leq \frac{ \nu 
 (\si,\tau)}{\al},$$
 holds with $p = 2 n - 2 + \al$.
 \end{cor}
  \section{A minimum principle for quasi-plurisubharmonic functions}

  Let $X$ a be compact K\"ahler manifold of dimension $n$ and $\w$ 
 a closed positive current on $X$ with bounded local potentials 
 such that the volume $Vol_{\w} (X) := \int_X \w^n > 0$.
 
 Let us fix  adenote by $d$ the geodesic metric on $X$ et denote by $h^p$ the Hausdorff
 content of dimension $p$ on the metric space $(X,d)$.
 
 Recall that a function $\vfi :X \lra  \R  \cup \{- 
 \infty\}$ is called an $\w-$plurisubharmonic function on $X$ if  
 $\vfi$ is upper semi-continuous on $X$ and satisfies the condition 
 $dd^c \vfi + \w \geq 0$ in the sense of currents on $X$. Therefore 
 each point of $X$ has a small neighbourhood $U \sub X$ which is 
 biholomorphic to a an euclidean ball in $\C^n$ so that $  
 v_U := \vfi + p_{U}$ is a \psh function on a neighbourhood of   
 $\ove U$,  where $p_U$ is a 
 local bounded potential for $\w$ on a neighbourhood of $\ove U$   
 i.e. $dd^c p_U = \w$ on a neighbourhood of $\ove U$.
 By compactness, $X$ can be covered by a finite number of such 
 domains. Let us  denote by $N = N (X)$ 
 the minimum number of such domains necessary to cover $X$.
 
 Observe that the Lelong numbers of an $\w-$\psh function $\vfi$ 
 are well defined by the formula 
 $\vt_{\vfi} (x) = \vt_v (x),$ where $v := \vfi + p$ is \psh on a 
 neighbourhood of $x$ and $p$ is a  bounded local potential of $\w$ 
 in a neighbourhood of $x$.
 
 Let us define the maximal Lelong number of $(X,\w)$ by the formula
 $$\vt (X,\w) := \sup \{\vt_{\vfi} (x) ; \vfi \in PSH (X,\w), x \in 
 X \}.$$
 It is clear that $\vt (X,\w) = \sup \{\vt_{\vfi} (x) ; \vfi \in 
 \mc P_0, x \in X \},$ where $\mc P_0 := \{\vfi \in PSH (X,\w);  
 \sup_X \vfi = 0\}.$
 Since $\mc P_0$ is compact fro the $L^1$ topology (see [17]),  
 a standard compactness argument shows that $\vt (X,\w) < + \infty$  (see [39]).
 
 It's well known that if $X$ is a projective manifold and $\w = 
 \w_{FS}$ is the Fubini-Study K\"ahler form on $X$ then 
 $\vt (X,\w) = \int_X \w^n$, which is equal to the degree of  
 algebraicity of $X$ (see [38]). We do not know whether this  
 formula is true in the general case.
 
 To state an analogue of Bernstein-Walsh inequalities which is the 
 couterpart of the classical Hadamard-three circle maximum  
 principle for the Lelong class $\mc L (\C^n)$, let us recall the  
 definition of the global capacity of Borel sets in $(X,\w)$ (see  
 [17]).
 
 For any Borel subset $K \sub X$, we define the  
 $\w-$capacity of $K$ in $X$ by
 $$T_{\w} (K) := exp (- \sup_X V_{K,\w}),$$
 where $V_{K,\w}$ is the extremal $\w-$\psh function associated to  
 $K$ and defined as follows:
 $$V_{K,\w} (x) := \sup \{\vfi (x) ; \vfi \in PSH (X,\w)\}, x \in 
 X.$$
 Then by [17], for any non pluripolar subset $K \sub X$, $T_{\w} (K) > 0$ is 
 the best constant such that 
 $$\vfi (x) \leq \sup_K \vfi - \log T_{\w} (K),  \forall \vfi \in PSH (X,\w), \forall x \in X.$$
 
 Now we can state the minimum principle for $\w-$\psh functions 
 analoguous the minimum principle for \psh functions. Namely we 
 obtain the following result.
 \begin{thm} For any $\nu > \vt (X,\w)$ there exists a constant $C 
 = C (X,\w,\nu) > 0$ and a real number $\eta_0 > 0$  small enough 
 such that for  any real number $0 < \al \leq 2$, any real number 
 $0 < \eta \leq \eta_0$ and any $\w-$\psh function $\vfi$ on $X$ 
 the following lower bound
  \begin{equation}
  \vfi (x) \geq \sup_{X} \vfi - \nu \log (C \slash \eta) 
  \label{eq:mp}
  \end{equation}
 holds for any $z \in X$ outside the union of a countable family  
 of balls of radii $(r_j)$ satisfying the following condition
 \begin{equation}
 \sum r_j^{2 n - 2 + \al} <   N \ 9^{n - 1} \eta^{\al} \slash \al.
 \label{eq:hc}
 \end{equation}
 \end{thm}
 \demo Take a covering of $X$ by $N = N (X)$ domains $(U_i)$ for 
 which there are domains $(U'i)$ biholomorphic to the unit ball of 
 $\C^n$ with $U_i \Sub U_i'\sub X$  for $1 \leq i \leq N$ and write 
 $v_i := \vfi_i + p_i$ where $p_i$ is a bounded local potential for 
 $\w$ on $U_i'$.
 
 We can apply the same method as in the proof of Theorem 6.1. Then from the fact that $p_i$ is a bounded \psh 
 function on $U'_i$, we obtain for each $i$ constants $a_i > 0,  
 b_i, c_i > 0$ such that the following lower bound 
 $$ \vfi (z) \geq  a_i \int_{\bd U_i'} \vfi d \si_{2 n - 1} - 
 \theta_\vfi (z,s) \log (3 \slash \eta) - b_i \ \mu_{\vfi} (U'_i) 
 \log (e \slash s) - c_i, $$
 holds for all $z \in U_i \sm E_i$, where $E_i \sub U_i$ is a Borel 
 set with $h_{2 n - 2 + \al}^{\eta} (E_i) < 9^{n - 1} \eta^{\al} 
 \slash \al$.
 By compactness of the $\mc P_0 := \{\vfi \in PSH (X,\w) ; \max_{X} 
 \vfi\} = 0$ (see [17]), it follows that $ \int_{\bd U_i'} \vfi d 
 \si_{j}$ and $ \mu_{\vfi} (U_i)$ are uniformly bounded for $\vfi  
 \in \mc P_0$. 
 Therefore taking $E := \bigcup_{1 \leq i \leq N} E_i$, we obtain  
 for a suitable constant $C > 0$ the following lower  bound
 $$\vfi  (z) \geq - \theta_{\vfi } (z,s) \log (3 \slash \eta) - C 
 \geq  - \nu \log (3 \slash \eta) - C, $$
 for any $z \in X \sm E$ and any $\vfi \in \mc P_0$.
 Since $h^{2 n - 2 + \al}_{\eta} (E) \leq \sum_i h_{2 n - 2 +  
 \al}^{\eta} (E_i) < N 9^{n - 1} \eta^{\al} \slash \al$ and $\vfi - 
 \max_X \vfi \in \mc P_0$ for any $\vfi \in PSH (X,\w)$, we obtain 
 the estimate of the theorem with an appropriate constant. \fin
 
 Now we want to prove the following comparison inequality.
 \begin{thm} For any $0 < \theta < \theta (X,\w),$ there exists 
 constant $C > 0$ such that for any $0 < \al \leq 2$ the following 
 comparison inequality
 $$ h^{2 n - 2 + \al} (K) \leq \frac{C}{\al} T_{\w} (K)^{\al 
 \theta},$$
 holds for ay Borel subset $K \sub X$.
 \end{thm}
 
 \noindent {\bf Problem:} Is the theorem still true with $\theta = 
 \theta (X,\w)$?
 \vskip 0.3 cm
 \noindent{\it Aknowlegments:} A lecture on a part of the  
 preliminary version of this paper was presented in jully 16, 2005  
 in Rabat at the conference held in honour to my colleague and friend Bensalem 
 Jennane. It's a pleasure for me to dedicate this work to him and 
 to thank him for his valuable collaboration and his kind  
 hospitality during all my numerous stays in 
 Rabat. I whould also like to thank Abdelhak Azhari, Ahmed Sebbar  
 and Alain Yger for interesting discussions on this subject and for  providing me interesting references.

 \vskip 0.3 cm  
 \noindent Universit\'e Paul Sabatier - Toulouse 3,  \\ 
 Laboratoire de Math\'ematiques Emile Picard,\\
 118 route de Narbonne, \\
 31062 Toulouse Cedex, France   \\ 
 e-mail: zeriahi $@$ math.ups-tlse.fr \\
 
 \end{document}